\documentclass[12pt]{amsart}
\usepackage{amsfonts}
\usepackage{amssymb}
\theoremstyle{definition}

\theoremstyle{remark}

\newcommand{\ds}{\displaystyle}

\begin{document}

Rendiconti di Matematica, Serie VII

Volume 18, Roma (1998), 151-166

\vspace{1.2in}

\centerline{\Large\bf Almost K\"ahler manifolds whose antiholomorphic}
\centerline{\Large\bf sectional curvature is pointwise constant
\footnote{ This paper has been partially supported by M.U.R.S.T. \\
Key Words and Phrases: {\it Almost K\"ahler manifold - Antiholomorphic sectional curvature 
- Complex space-form }\\ 
 A.M.S. Classification: 53C15 - 53C55 - 53B21}}

\vspace{0.6in}
\centerline{\large\bf M. FALCITELLI -- A. FARINOLA -- O.T. KASSABOV}

\vspace{0.6in}
{\sl Abstract: We prove that an almost K\"ahler manifold $(M,g,J)$ with ${\rm dim}\,M\ge8$
and pointwise constant antiholomorphic sectional curvature is a complex space-form. }

\vspace{0.4in}
{\bf 1 -- Introduction and preliminaries} 

\vspace{0.1in}
Let $(M,g,J)$ be a $2n$-dimensional almost Hermitian manifold. A 2-plane $\alpha$ in the
tangent space $T_xM$ at a point $x$ of $M$ is antiholomorphic if it is orthogonal to $J\alpha$.

The manifold $(M,g,J)$ has pointwise constant antiholomorphic sectional curva\-ture (p.c.a.s.c.)
$\nu$ if, at any point $x$, the Riemannian sectional curvature $\nu(x)=K_x(\alpha)$ is
independent on the choice of the antiholomorphic 2-plane $\alpha$ in $T_xM$.

If $(g,J)$ is a K\"ahler structure, the previous condition means that $(M,g,J)$ is a complex
space-form, i.e. a K\"ahler manifold with constant holomorphic sectional curvature $\mu=4\nu$ ([2]).
Moreover, the Riemannian curvature tensor $R$ satisfies:
\vspace{0.1in}
$$
	R=\nu(\pi_1+\pi_2)\,,    \leqno (1.1)
$$
\vspace{0.1in}
$\nu$ being  a constant function and $\pi_1,\,\pi_2$ the tensor fields such that:
$$
	\begin{array}{c}\vspace{0.1in}
		\pi_1(X,Y,Z,W)=g(X,Z)g(Y,W)-g(Y,Z)g(X,W) \,;  \\ \vspace{0.1in}
		\pi_2(X,Y,Z,W)=2g(JX,Y)g(JZ,W)+g(JX,Z)g(JY,W)   \\ \vspace{0.1in}
		-g(JY,Z)g(JX,W)\,.
	\end{array}    \leqno(1.2)
$$

According to [16], for any (0,2)-tensor field $S$, we consider the (0,4)-tensor fields
$\phi(S)$, $\psi(S)$ defined by:
$$
	\begin{array}{r}\vspace{0.1in}
		\phi(S)(X,Y,Z,W)=g(X,Z)S(Y,W)+g(Y,W)S(X,Z)  \\ \vspace{0.1in}
										-g(X,W)S(Y,Z)-g(Y,Z)S(X,W) \,,  \\ \vspace{0.1in}
		\psi(S)(X,Y,Z,W)=2g(X,JY)S(Z,JW)+2g(Z,JW)S(X,JY)   \\ \vspace{0.1in}
		                +g(X,JZ)S(Y,JW)+g(Y,JW)S(X,JZ)    \\ \vspace{0.1in}
		                -g(X,JW)S(Y,JZ)-g(Y,JZ)S(X,JW) \,.
\end{array}    \leqno(1.3)
$$
A generalization of (1.1) is obtained by G. Ganchev ([5]). In fact, he proves that
the almost Hermitian manifold $(M,g,J)$ has p.c.a.s.c. $\nu$ iff

$$
	R=\frac{1}{2(n+1)}\psi(\rho^*(R))+\nu\pi_1-\frac{2(n+1)\nu+\tau^*(R)}{2(n+1)(2n+1)}\pi_2  \leqno(1.4)
$$

\vspace{0.1in}\noindent
$\rho^*(R),\,\tau^*(R)$ respectively denoting the $*$-Ricci tensor and the $*$-scalar 
curvature.

The previous formula allows to relate the symmetric part of $\rho^*(R)$ to the Ricci tensor 
$\rho(R)$ and thus $\tau^*(R)$ to the scalar curvature $\tau(R)$.

Indeed, putting

$$
	L_3R(X,Y,Z,W)=R(JX,JY,JZ,JW)  \,,   \leqno (1.5)
$$

\vspace{0.1in}\noindent
one has: 

$$
	\rho^*(R+L_3R)(X,Y)=\rho^*(R)(X,Y)+\rho^*(R)(Y,X)=\rho^*(R+L_3(R))(JX,JY) \,,
$$

\vspace{0.1in}\noindent
and (1.4) implies:

$$
	\rho^*(R+L_3R)=\frac23(n+1)\rho(R)-\frac{(n+1)\tau(R)-3\tau^*(R)}{3n}g \,,  \leqno(1.6)
$$

$$
	8n(n^2-1)\nu=(2n+1)\tau(R)-3\tau^*(R) \,.   \leqno (1.7)
$$

\vspace{0.1in}
Another characterization of the p.c.a.s.c. condition can be obtained regarding the 
Riemannian curvature tensor as a section of the vector bundle $\mathcal R(M)$ of 
the algebraic curvature tensor fields on $M$. According to the splitting
$\mathcal R(M)=\oplus_{1\le i\le10}  \mathcal W_i(M)$ considered in [16], the formula
(1.4) can be interpreted in terms of the vanishing of suitable $\mathcal W_i$-projections
$p_i(R)$ of $R$.

More precisely, an application of the Theorem 8.1 in [16] yields to the following result.

\vspace{0.3in}
P\textsc{roposition} 1.1. {\it Let $(M,J,g)$ be an almost Hermitian manidold. If ${\rm dim}\,M=4$,
$(M,g,J)$ has p.c.a.s.c. iff $p_3(R)=p_7(R)=p_8(R)=0$.
If ${\rm dim}\,M\ge 6$, then $(M,g,J)$ has p.c.a.s.c. iff $p_3(R)=p_6(R)=p_7(R)=p_8(R)=p_{10}(R)=0$
and (1.6) holds.}

\vspace{0.3in}
Combining with the Theorem 18 in [4], one has:

\vspace{0.3in}
P\textsc{roposition} 1.2. {\it Let $(M,J,g)$ be an almost Hermitian manidold with p.c.a.s.c.
Then $g$ is an Einstein metric iff $(M,g,J)$ has pointwise constant holomorphic 
sectional curvature.} 

\vspace{0.3in}
The classification of the almost Hermitian manifolds with p.c.a.s.c. is still an open 
problem, even if nowadays several partial results are known.

In [1] V. Apostolov, G, Ganchev and S. Ivanov classify the compact Hermitian surfaces 
with constant antiholomorphic curvature. Moreover, they construct an example of conformal 
K\"ahler surface with p.c.a.s.c. $\nu$, the function $\nu$ being non-constant. Thus,
the Schur's lemma of antiholomorphic type is not valid in the 4-dimensional case.

Furthermore, the third autor of the present paper has already solved the above-mensioned
problem for $2n$-dimensional, $n\ge 3$, connected, $\mathcal R_3-$manifolds, i.e. almost 
Hermitian manifolds such that $R=L_3R$ (equivalently, $p_8(R)=p_9(R)=p_{10}(R)=0$).

In fact, any connected $\mathcal R_3$-manifold $M$ with p.c.a.s.c. and ${\rm dim}\,M\ge 6$
has constant antiholomorphic sectional curvature ([9]) and turns out to be a real
space-form or a complex space-form ([10]).

This result alows the classification of nearly K\"ahler as well as locally conformal K\"ahler
manifolds with p.c.a.s.c. In fact, any nearly K\"ahler manifold is a $\mathcal R_3$-manifold
([7]). Since for a locally conformal K\"ahler manifolds the projections $p_9(R)$ vanishes, 
the locally conformal K\"ahler manifolds with p.c.a.s.c. turn out to be $\mathcal R_3$-manifolds ([3]).

Moreover, combining the results stated in [10] and [13], any connected $\mathcal R_3$-almost 
K\"ahler manifold with p.c.a.s.c. and ${\rm dim}\,M\ge 6$ turns out to be a complex
space-form.

Since the projection $p_9(R)$, a priori, does not vanish in the almost K\"ahler case
the classification of the almost K\"ahler manifolds with p.c.a.s.c is meaningful.

We recall the almost K\"ahler condition, i.e.:

$$
	\underset{V,X,Y}\sigma(\nabla_V\omega)(X,Y)=0  \,,  \leqno (1.8)
$$

\vspace{0.1in}\noindent
$\sigma$ denoting the cyclic sum and $\nabla\omega$ the covariant derivative of the
fundamental 2-form $\omega \ $ ($\omega(X,Y)=g(JX,Y)$) with respect to the Levi-Civita
connection $\nabla$.

Moreover, (1.8) implies:

$$
	(\nabla_XJ)Y+(\nabla_{JX}J)JY=0  \,;   \leqno (1.9)
$$

$$
	\sum_i (\nabla_{e_i}J)e_i=0  \,,   \leqno (1.10)
$$

\vspace{0.1in}\noindent
for any local orthonormal frame $\{ e_i \}_{1\le i\le 2n}$.

In this paper we state  the following theorem, whose proof is divided into several steps.

\vspace{0.3in}
T\textsc{heorem} 1. {\it Let $(M,g,J)$ be a $2n$-dimensional, $n\ge 4$, connected,
almost K\"ahler manifold. If $(M,g,J)$ has pointwise constant antiholomorphic
sectional curvature, then $(M,g,J)$ is a complex space-form.}

\vspace{0.5in}
{\bf 2 -- Some auxiliary lemmas} 

\vspace{0.1in}
Given a $2n$-dimensional almost Hermitian manifold $(M,g,J)$, the tensor field

$$
	Q=\frac16 \rho(R)+\frac1{4(n+1)} \rho^*(R-L_3R)   \leqno (2.1)
$$

\vspace{0.1in}\noindent
is, in general, neither symmetric nor skew-symmetric, since $\rho(R)$,
$\rho^*(R-L_3R)$ respectively determine its symmetric, skew-symmetric components.
Moreover, we assume that $(M,g,J)$ has p.c.a.s.c.; then the formula (1.6)
implies:

$$
	Q(JX,JY)=Q(Y,X)  \,,   \leqno (2.2)
$$

\vspace{0.1in}\noindent
and thus one has:

$$
	Q((\nabla_VJ)X,JY)=Q(Y,(\nabla_VJ)JX)) \,;   \leqno (2.3)
$$

$$
	\begin{array}{r}  \vspace{0.1cm}
		(\nabla_VQ)(JX,JY)=(\nabla_VQ)(Y,X)-Q((\nabla_VJ)X,JY)  \\ \vspace{0.1cm}
					-Q(JX,(\nabla_VJ)Y) \,;
	\end{array}   \leqno (2.4)
$$

$$
	\sum_i Q((\nabla_VJ)e_i,Je_i)=-\sum_iQ(Je_i,(\nabla_VJ)e_i,)   \,,  \leqno (2.5)
$$

\vspace{0.1in}\noindent
for any local orthonormal frame $\{ e_i \}_{1\le i\le 2n}$.

If $(g,J)$ is an almost K\"ahler structure, (1.8) and (1.9) imply also

$$
	\sum_iQ(V,e_i)(\nabla_{e_i}\omega)(Y,X)=Q(V,(\nabla_XJ)Y-(\nabla_YJ)X)  \,;   \leqno(2.6)
$$

$$
	2\sum_i Q(Je_e,(\nabla_{e_i}J)V)=\sum_iQ(Je_i,(\nabla_VJ)e_i) \,.   \leqno (2.7)
$$

\vspace{0.1in}
Now we observe that (2.1), (1.6) and (1.7) allow to rewrite (1.4) as follows:

$$
	R=\psi(Q)+\nu\pi_1-\frac{2n-1}3\nu\pi_2  \,. \leqno (2.8)
$$

\vspace{0.1in}
By means of (2.8) and the second Bianchi identity, we will state some properties
of $Q$ and $\nabla Q$ useful for the proof of the Theorem 1.

First of all, from (2.8), one has:

$$
	\begin{array} {r} \vspace{0.1in}
		(\nabla_VR)(X,Y,Z,W)=2g(X,JY) \{ (\nabla_VQ)(Z,JW)  \\ \vspace{0.1in}
		+Q(Z,(\nabla_VJ)W) \} +2g(Z,JW) \{ (\nabla_VQ)(X,JY)  \\ \vspace{0.1in}
		+Q(X,(\nabla_VJ)Y) \} + g(X,JZ) \{ (\nabla_VQ)(Y,JW)  \\ \vspace{0.1in}
		+Q(Y,(\nabla_VJ)W) \} + g(Y,JW) \{ (\nabla_VQ)(X,JZ)  \\ \vspace{0.1in}
		+Q(X,(\nabla_VJ)Z) \} - g(Y,JZ) \{ (\nabla_VQ)(X,JW)  \\ \vspace{0.1in}
		+Q(X,(\nabla_VJ)W) \} - g(X,JW) \{ (\nabla_VQ)(Y,JZ)  \\ \vspace{0.1in}
		+Q(Y,(\nabla_VJ)Z) \} +2(\nabla_V\omega)(Y,X)Q(Z,JW)  \\ \vspace{0.1in}
	\end{array}  \leqno (2.9)
$$ 		

$$
	\begin{array} {r} \vspace{0.1in}		
		+2(\nabla_V\omega)(W,Z)Q(X,JY) + (\nabla_V\omega)(Z,X)Q(Y,JW) \\ \vspace{0.1in}
		+ (\nabla_V\omega)(W,Y)Q(X,JZ) - (\nabla_V\omega)(Z,Y)Q(X,JW) \\ \vspace{0.1in}
		- (\nabla_V\omega)(W,X)Q(Y,JZ)  \\ \vspace{0.1in}
		\ds +V(\nu)\big(\pi_1-\frac{2n-1}3 \pi_2\big)(X,Y,Z,W)  \\ \vspace{0.1in}
		\ds -\frac{2n-1}3\nu \{ 2g(X,JY)(\nabla_V\omega)(W,Z)  \\ \vspace{0.1in}
		+2g(Z,JW)(\nabla_V\omega)(Y,X) + g(X,JZ)(\nabla_V\omega)(W,Y)  \\ \vspace{0.1in}
		+ g(Y,JW)(\nabla_V\omega)(Z,X) - g(X,JW)(\nabla_V\omega)(Z,Y) \\ \vspace{0.1in}
		-g(Y,JZ)(\nabla_V\omega)(W,X) \} \,.
	\end{array} 
$$ 

\vspace{0.3in}
L\textsc{emma} 2.1. {\it Let $(M,g,J)$ be a $2n$-dimensional $(n\ge 2)$ almost-K\"ahler
manifold with p.c.a.s.c. The covariant derivative $\nabla Q$ is given by:

$$
	\begin{array} {r} \vspace{0.1in}
		2(n+1)(2n-1)(\nabla_VQ)(X,JY)=(2n+3)(Q(Y,(\nabla_XJ)V)  \\ \vspace{0.1in}
		-Q(X,(\nabla_YJ)V)+ (4n+3)Q(V,(\nabla_XJ)Y-(\nabla_YJ)X)  \\ \vspace{0.1in}
		-Q(Y,(\nabla_VJ)X) -(4n^2+2n-3)Q(X,(\nabla_VJ)Y) \\ \vspace{0.1in}
		\ds +g(X,JY)\big\{2n \sum_i Q(Je_i,(\nabla_VJ)e_i)  \\ \vspace{0.1in}
		\ds + \frac43(n+1)(n-2)V(\nu) +\frac{2n-1}6 V(\tau(R)) \big\} \\ \vspace{0.1in}
		\ds +g(X,JV)\big\{\frac{4n-1}2 \sum_i Q(Je_i,(\nabla_YJ)e_i)  \\ \vspace{0.1in}
		\ds - \frac23(n+1)(2n^2-4n+3)Y(\nu) +\frac{2n-1}6 Y(\tau(R)) \big\} \\ \vspace{0.1in}
		\ds -g(Y,JV)\big\{\frac{4n-1}2 \sum_i Q(Je_i,(\nabla_XJ)e_i)  \\ \vspace{0.1in}
		\ds - \frac23(n+1)(2n^2-4n+3)X(\nu) +\frac{2n-1}6 X(\tau(R)) \big\} \\ \vspace{0.1in}
		-2(n+1)\{ JX(\nu)g(Y,V)-JY(\nu)g(X,V) \}   \\ \vspace{0.1in}
		+\ds\frac13 (n+1)(\tau(R)-2(2n-1)^2\nu)(\nabla_V\omega)(X,Y) \,,
	\end{array}  \leqno (2.10)
$$

\vspace{0.1in}\noindent
where $\{ e_i \}_{1\le i\le 2n}$ is a local orthonormal frame.}

\vspace{0.2in}
P\textsc{roof}. In fact, by the second Bianchi identity, we have:
$$
	\underset{(V,X,Y)}\sigma\sum_i (\nabla_VR)(X,Y,e_i,Je_i)=0 \,,
$$
which, combined with (2.9), (2.4) and (1.8), yields to:

$$
	\begin{array} {r}\vspace{0.1in}
		\ds 2(n+1)\underset{(V,X,Y)}\sigma(\nabla_VQ)(X,JY)=\underset{(V,X,Y)}\sigma \big\{ Q(Y,(\nabla_VJ)X) \\ \vspace{0.1in}
		-\ds(2n+3)Q(X,(\nabla_VJ)Y)+g(JX,Y) \big[ \frac16V(\tau(R))  \\ \vspace{0.1in}
		\ds -\frac43(n^2-1)V(\nu) + \sum_i Q(Je_i,(\nabla_VJ)e_i) \big]\big\}  \ .
	\end{array}  \leqno (2.11)
$$

\vspace{0.1in}
Moreover, by the second Bianchi identity, we obtain:

$$
	\sum_{i,q} \{ 2(\nabla_{e_i}R)(V,e_q,Je_i,Je_q) - (\nabla_VR)(e_i,e_q,Je_i,Je_q) \} =0 \,,
$$

\vspace{0.1in}\noindent
which, combined with (2.9), (2.5) and (2.7) implies:

$$
	\begin{array} {r}\vspace{0.1in}
		\ds \sum_i(\nabla_{e_i}Q)(V,e_i)=\frac{4n+1}{4(n+1)} \sum_i Q(Je_i,(\nabla_VJ)e_i) \\
		       \ds +\frac{n}{6(n+1)}V(\tau(R))-\frac23(n-1)^2V(\nu) \,.
	\end{array}  \leqno (2.12)
$$

\vspace{0.1in}
This formula, with (2.4), (2.6), (1.10) and the condition

$$
	\sum_i \{ (\nabla_VR)(e_i,Le_i,X,Y)+2(\nabla_{e_i}R)(Je_i,V,X,Y) \} =0
$$

\vspace{0.1in}\noindent
yields to

$$
	\begin{array} {r}\vspace{0.2cm}
		2n(\nabla_VQ)(X,JY)+(\nabla_XQ)(Y,JV)+(\nabla_YQ)(V,JX)=2Q(V,(\nabla_XJ)Y)  \\ \vspace{0.1in}
		-3Q(V,(\nabla_YJ)X)-2nQ(X,(\nabla_VJ)Y)-Q(X,(\nabla_YJ)V)  \\ \vspace{0.1in}
		\ds +g(X,JY) \big\{ \frac{2n-1}{2(n+1)} \sum_iQ(Je_i,(\nabla_VJ)e_i)+\frac23 (2n-3)V(\nu) \\ \vspace{0.1in}
		\ds +\frac{n-1}{6(n+1)}V(\tau(R))\big\} +g(V,JY) \big\{ \frac{4n+1}{4(n+1)} \sum_i Q(Je_i,(\nabla_XJ)e_i) \\ \vspace{0.1in}
		\ds -\frac13 (2n^2-2n+1)X(\nu) + \frac{n}{6(n+1)}X(\tau(R)) \big\}  \\ \vspace{0.1in}
		\ds-g(V,JX) \big\{ \frac{4n+1}{4(n+1)} \sum_i Q(Je_i,(\nabla_YJ)e_i) \\ \vspace{0.1in}
		\ds -\frac13 (2n^2-2n+1)Y(\nu) + \frac{n}{6(n+1)}Y(\tau(R)) \big\}  \\ \vspace{0.1in}
		-JX(\nu)g(Y,V)+JY(\nu)g(X,V)  \\ \vspace{0.1in}
		+\ds \frac16(\tau(R)-2(2n-1)^2\nu)(\nabla_V\omega)(X,Y) \,.
	\end{array}
$$

\vspace{0.1in}
Thus, combining with (2.11), one proves the statement.

\vspace{0.3in}
L\textsc{emma} 2.2. {\it In the hypothesis of the Lemma 2.1, when $n \ne 3$, one has:}
$$
	\sum_i Q((Je_i,(\nabla_VJ)e_i) = \frac43(n^2-1)V(\nu)  \,.  \leqno (2.13)
$$

\vspace{0.2in}
P\textsc{roof}. In fact, the second Bianchi identity and (2.1) give:

$$
	V(\tau(R))=2\sum_i (\nabla_{e_e}\rho(R))(V,e_i) =
	  6\sum_i\{(\nabla_{e_i}Q)(V,e_i)+(\nabla_{e_i}Q)(e_i,V)\} \,.
$$

\vspace{0.1in}\noindent
Moreover, the formulas (2.10), (2.12), (2.7), (1.8), (1.9), (1.10) imply:

$$
	\sum_i\{(\nabla_{e_i}Q)(V,e_i)+(\nabla_{e_i}Q)(e_i,V)\}=\frac16V(\tau(R))
$$

$$
	+\frac{n-2}{2n-1} \big\{ \sum_i Q(Je_i,(\nabla_VJ)e_i)-\frac43(n^2-1)V(\nu) \big\} \,,
$$

\vspace{0.1in}\noindent
and then the statement.

\vspace{0.3in}
P\textsc{roposition} 2.1. {\it In the hypothesis of the Lemma 2.1, if $n\ge 3$, one has:}

$$
	\begin{array}{r} \vspace{0.1in}
		4(2n-3)(Q(X,(\nabla_YJ)W)-Q(Y,(\nabla_XJ)W) \\ \vspace{0.1in}
		+4n(Q(X,(\nabla_WJ)Y) - Q(Y,(\nabla_WJ)X)) \\ \vspace{0.1in}
		-4(n-3)Q(W,(\nabla_XJ)Y-(\nabla_YJ)X) +\tau(R)(\nabla_W\omega)(X,Y)  \\ \vspace{0.1in}
		\ds -\frac83(2n^2-4n+3)(X(\nu)g(JY,W)  \\ \vspace{0.1in}
		-Y(\nu)g(JX,W)+2W(\nu)g(X,JY)  \\ \vspace{0.1in}
		+8n(n-2)(JX(\nu)g(Y,W)  \\ \vspace{0.1in}
		-JY(\nu)g(X,W))=0 \,. 
	\end{array}  \leqno (2.14)
$$

\vspace{0.2in}
P\textsc{roof.} The Lemmas 2.1 and 2.2, the formula (2.4) and the condition:

$$
	\sum_i \{(\nabla_{e_i}R)(X,Y,e_i,W)+(\nabla_YR)(e_i,X,e_i,W)-(\nabla_XR)(e_i,Y,e_i,W)\}=0
$$

imply the vanishing of the tensor field $S$ defined by:

$$
	\begin{array}{r} \vspace{0.1in}
		S(W,X,Y)=4(2n^2-3) \{Q(X,(\nabla_YJ)W)-Q(Y,(\nabla_XJ)W) \} \\ \vspace{0.1in}
		-4n \{ Q((\nabla_YJ)W,X) - Q((\nabla_XJ)W,Y) \} \\ \vspace{0.1in}
		+2(2n^2+3n+3)\{  Q(X,(\nabla_WJ)Y) - Q(Y,(\nabla_WJ)X) \} \\ \vspace{0.1in}
		-2(n+3) \{ Q((\nabla_WJ)Y,X)-Q((\nabla_WJ)X,Y) \} \\ \vspace{0.1in}
		+2(n-3) \{ Q(W,(\nabla_XJ)Y-(\nabla_YJ)X) -(2n+3)Q(\nabla_XJ)Y \\ \vspace{0.1in}
		 -(\nabla_YJ)X,W)\}+(n+1)\tau(R) (\nabla_W\omega)(X,Y) \\ \vspace{0.1in}
		-4(n+1)(2n^2-2n-3)\{ X(\nu)g(JY,W)-Y(\nu)g(JX,W)\}  \\ \vspace{0.1in}
		\ds +\frac 43(n+1)(4n^2-4n+3)\{ JX(\nu)g(Y,W)  \\ \vspace{0.1in}
		-JY(\nu)g(X,W)-2W(\nu)g(X,JY) \} \,.
	\end{array}  
$$

In particular, by means of (1.9) and (2.3), the conditions:

$$
	S(W,X,Y)-S(W,JX,JY)-S(JW,JX,Y)-S(JW,X,JY)=0 \,;
$$

$$
	S(W,X,Y)-S(W,JX,JY)+S(JW,JX,Y)+S(JW,X,JY)=0 \,,
$$

turn out to be equivalent to:

$$
	\begin{array}{r} \vspace{0.1in}
		2(n-3)\{ Q(W,(\nabla_XJ)Y-(\nabla_YJ)X) + Q((\nabla_XJ)Y  \\ \vspace{0.1in}
		-(\nabla_YJ)X,W)\} -2(2n-3) \{ Q(X,(\nabla_YJ)W)  \\ \vspace{0.1in}
		+Q((\nabla_YJ)W,X)-Q(Y,(\nabla_XJ)W)-Q((\nabla_XJ)W,Y) \} \\ \vspace{0.1in}
		-2n \{ Q(X,(\nabla_WJ)Y) +Q((\nabla_WJ)Y,X)   \\ \vspace{0.1in}
		-Q(Y,(\nabla_WJ)X)-Q((\nabla_WJ)X,Y)  \\ \vspace{0.1in}
		-\tau(R)(\nabla_W\omega)(X,Y) = 0 \,;
	\end{array}  \leqno (2.15)
$$

$$
	\begin{array}{r} 
		(n-3)\big\{ Q(W,(\nabla_YJ)X-(\nabla_XJ)Y) - Q((\nabla_YJ)X  \\ \vspace{0.1in}
		-(\nabla_XJ)Y,W)  +\ds\frac23(n+1)[X(\nu)g(JY,W)  \\ \vspace{0.1in}
		- Y(\nu)g(JX,W) +JX(\nu)g(Y,W)-JY(\nu)g(X,W) ]\big\} =0 \,. \vspace{0.1in}
	\end{array}  \leqno (2.16)
$$

Thus, if $n=3$, the statement follows from (2.15) combined with the
condition $S=0$. If $n>3$, (2.16) implies also, with suitable change of
the involved variables, the relation:

$$
	\begin{array}{r} \vspace{0.1in}
		Q((\nabla_WJ)X,Y) - Q((\nabla_WJ)Y,X) = Q((\nabla_XJ)W,Y) \\ \vspace{0.1in}
		-Q((\nabla_YJ)W,X) + Q(Y,(\nabla_WJ)X-(\nabla_XJ)W) \\ \vspace{0.1in}
		-Q(X,(\nabla_WJ)Y-(\nabla_YJ)W) +\ds \frac23(n+1)\big\{ X(\nu)g(Y,JW) \\  \vspace{0.1in}
		- Y(\nu)g(X,JW) + 2 W(\nu)g(X,JY)  \\ \vspace{0.1in}
		+JX(\nu)g(Y,W) -JY(\nu)g(X, W) \big\} \,.
	\end{array}  \leqno (2.17)
$$

Thus, applying (2.17) and (2.16), the relation (2.15) yields to:

$$
	\begin{array}{r} \vspace{0.1in}
		6(n-1)\{ Q((\nabla_XJ)W,Y) - Q((\nabla_YJ)W,X) \}  \\ \vspace{0.1in}
		=2(n-3) \{ Q(X,(\nabla_YJ)W) - Q(Y,(\nabla_XJ)W) \} \\ \vspace{0.1in}
		+4n \{ Q(X,(\nabla_WJ)Y) - Q(Y,(\nabla_WJ)X) \}  \\ \vspace{0.1in}
		-4(n-3) Q(W,(\nabla_XJ)Y-(\nabla_YJ)X)  \\ \vspace{0.1in}
		+\tau(R)(\nabla_W\omega)(X,Y) +\ds \frac43(n+1)\{ (2n-3)X(\nu)g(JY,W) \\ \vspace{0.1in}
		-Y(\nu)g(JX,W) -2n W(\nu)g(X,JY) \} \\ \vspace{0.1in}
		-4(n+1) \{ JX(\nu)g(Y,W) - JY(\nu)g(X,W) \} \,. \vspace{0.1in}
	\end{array}  \leqno (2.18)
$$

Moreover, via (2.17), (2.18) and (2.16), with a direct computation, one has:

$$
	\begin{array}{rl} \vspace{0.1in}
		S(W,X,Y) = & \ds \frac{n(n+1)}{n-1} \big\{ 4(2n-3) ( Q(X,(\nabla_YJ)W) \\ \vspace{0.1in}
		           &- Q(Y,(\nabla_XJ)W) +4n( Q(X,(\nabla_WJ)Y) \\ \vspace{0.1in}
	          	 &-Q(Y,(\nabla_WJ)X) -4(n-3)Q(W,(\nabla_XJ)Y-(\nabla_YJ)X) \\ \vspace{0.1in}
		           &+\tau(R)(\nabla_W\omega)(X,Y) +8n(n-2)(JX(\nu)g(Y,W)  \\ \vspace{0.1in}
		           &-JY(\nu)g(X,W) +\ds\frac83(2n^2-4n+3)( X(\nu)g(Y,JW) \\ \vspace{0.1in}
		           &-Y(\nu)g(X,JW) - 2 W(\nu)g(X,JY) ) \big\}  \,.
	\end{array}  
$$

Therefore, the vanishing of $S$ implies the statement.

\vspace{0.3in}
P\textsc{roposition} 2.2. {\it In the hypothesis of the Lemma 2.1, if $n\ge 4$, one has: }

$$
	\begin{array}{rl} \vspace{0.1in}
		Q(X,(\nabla_YJ)V) &-\ds Q((\nabla_YJ)V,X)=\frac23 \{(2n-1)( Y(\nu)g(JV,X) \\ \vspace{0.1in}
		                  &+ JY(\nu)g(V,X)) +(n-2)(V(\nu)g(JY,X) \\ \vspace{0.1in}
	          	        &+ JV(\nu)g(Y,X)) \} \,. \\ 
	\end{array}  \leqno (2.19)
$$

P\textsc{roof.} We consider the (0,3)-tensor field $T$ such that:

$$
	T(V,X,Y)=Q(V,(\nabla_XJ)Y-(\nabla_YJ)X) \,.
$$

Since $T$ satisfies:
$$
	T(V,X,Y)=- T(V,Y,X)=-T(V,JX,JY) \,,
$$

\noindent
$T$ can be regarded as a section of the vector bundle $\mathcal W(M)$
whose fibre, at any point $x$ of $M$, is the linear space $\mathcal W_x$
considered in [8].

According to the splitting $\ds\mathcal W(M)=\underset{1\le i\le4}   \oplus \ \mathcal W_i(M)$
defined in [8], we define by $q_1(T)$ the $\mathcal W_1$-projection of $T$;
it is the skew-symmetric tensor field such that:
$$
	6q_1(T)(V,X,Y)=\underset{(V,X,Y)}\sigma \ (T(V,X,Y)-T(JV,JX,Y)) \,.
$$

Since $n\ge 4$, applying (2.16) and then (2.14), one obtains:

$$
	\begin{array}{rl} \vspace{0.1in}
		3q_1(T)(V,X,Y) = & \underset{(V,X,Y)}\sigma \ ( Q(V,(\nabla_XJ)Y-(\nabla_YJ)X)    \\ \vspace{0.1in}
		                  &+ \ds \frac23 (n+1)V(\nu)g(X,JY)) \\ \vspace{0.1in}
	          	       = &\ds\frac1n \big\{ 3Q(V,(\nabla_XJ)Y-(\nabla_YJ)X) \\ \vspace{0.1in}
	          	        & +3(n-1)( Q(X,(\nabla_YJ)V)  \\ \vspace{0.1in}
	          	        & -Q(Y,(\nabla_XJ)V) +\ds\frac14 \tau(R)(\nabla_V\omega)(X,Y)  \\ \vspace{0.1in}
	          	        & +2n(n-2)(JX(\nu)g(V,Y)-JY(\nu)g(V,X)) \\ \vspace{0.1in}
	          	        & +2(n^2-n+1)(X(\nu )g(JV,Y)-Y(\nu)g(JV,X)) \\ \vspace{0.1in}
	          	        & -2(n-1)(n-2)V(\nu)g(X,JY) \big\} \,.
	\end{array}  
$$

Then, the condition: $q_1(T)(V,X,Y)+q_1(T)(JY,JX,Y)=0$ combined with (1.9), (2.3), (2.16)
proves the statement.

\vspace{0.5in}
{\bf 3 -- The proof of the Theorem 1} 

\vspace{0.1in}
To the Riemannian curvature of a manifold satisfying the hypothesis of the Theorem 1,
we apply the second Bianchi identity in the form:

$$
	\begin{array}{c} \vspace{0.1in}
			\underset{(V,X,Y)}\sigma \ \{(\nabla_VR)(X,Y,Z,W) + (\nabla_VR)(X,Y,JZ,JW) \}  \\ \vspace{0.1in}
			+\underset{(V,JX,JY)}\sigma  \{ (\nabla_VR)(JX,JY,Z,W)   \\ \vspace{0.1in}
			+(\nabla_VR)(JX,JY,JZ,JU) \} =0  \,.
	\end{array} \leqno (3.1)
$$

The complete expression of the first member in (3.1), evaluated by means of (2.9),
is a tensor field which contains four blocks of terms, respectively depending on
$ g \otimes(\nabla Q+Q(.,\nabla J))$, $\nabla\omega \otimes Q$, 
$d\nu \otimes(\pi_1 -\frac{2n-1}3\pi_2)$, $g\otimes\nabla\omega$.

Since $(g,J)$ is an almost K\'ahler structure, the whole term in $g\otimes \nabla\omega$
vanishes, while only the skew-symmetric component of $Q$, i.e. $\rho^*(R-L_3R)$, is
involved in the block depending on $\nabla\omega\otimes Q$.

After a quite long computation, applying the Lemmas 2.1 and 2.2 and then the 
Proposition 2.2, the whole expression in $ g \otimes(\nabla Q+Q(.,\nabla J))$
turns out to depend only on $d\nu \otimes g \otimes g$. Thus, the condition (3.1) 
is equivalent to
$$
	\begin{array}{l}  \vspace{0.1in}
		\ds\frac3{4(n+1)} \{\rho^*(R-L_3R)(X,Z)(\nabla_W\omega)(JY,V)   \\ \vspace{0.1in}
		\qquad -\rho^*(R-L_3R)(Y,Z)(\nabla_W\omega)(JX,V)   \\ \vspace{0.1in}
		\qquad -\rho^*(R-L_3R)(JX,Z)(\nabla_W\omega)(Y,V)   \\ \vspace{0.1in}
		\qquad +\rho^*(R-L_3R)(JY,Z)(\nabla_W\omega)(X,V) -\rho^*(R-L_3R)(X,W)\times  \\ \vspace{0.1in}
		\qquad \times(\nabla_Z\omega)(JY,V) + \rho^*(R-L_3R)(Y,W)(\nabla_Z\omega)(JX,V)   \\ \vspace{0.1in}
		\qquad +\rho^*(R-L_3R)(JX,W)(\nabla_Z\omega)(Y,V)   -\rho^*(R-L_3R)(JY,W)\times   \\ \vspace{0.1in}
		\qquad \times (\nabla_Z\omega)(X,V) \} -X(\nu)\{ \pi_1(V,Y,Z,W) + \pi_1(V,Y,JZ,JW) \\ \vspace{0.1in}
		\qquad +2g(Y,JV)g(Z,JW) \} +Y(\nu) \{ \pi_1(V,X,Z,W)  \\ \vspace{0.1in}
		\qquad +\pi_1(V,X,JZ,JW) +2g(X,JV)g(Z,JW) \}  \\ \vspace{0.1in}
		\qquad -JX(\nu) \{ \pi_1(V,JY,Z,W)   \\ \vspace{0.1in}
		\qquad -\pi_1(JV,Y,Z,W)+2g(Y,V)g(Z,JW) \}  \\ \vspace{0.1in}
		\qquad +JY(\nu) \{ \pi_1(V,JX,Z,W)   \\ \vspace{0.1in}
		\qquad -\pi_1(JV,X,Z,W)+2g(X,V)g(Z,JW) \}  \\ \vspace{0.1in}
		\qquad +2V(\nu) \{ \pi_1(X,Y,Z,W)   \\ \vspace{0.1in}
		\qquad +\pi_1( X,Y,JZ,JW)-2g(X,JY)g(Z,JW) \}  \\ \vspace{0.1in}
		\qquad +2W(\nu) \{ \pi_1(X,Y,Z,V)   \\ \vspace{0.1in}
		\qquad +\pi_1( X,Y,JZ,JV)-2g(X,JY)g(Z,JV) \}  \\ \vspace{0.1in}
		\qquad -2Z(\nu) \{ \pi_1(X,Y,W,V)   \\ \vspace{0.1in}
		\qquad +\pi_1( X,Y,JW,JV)-2g(X,JY)g(W,JV) \}  \\ \vspace{0.1in}
		\qquad -2JW(\nu) \{ \pi_1(X,Y,Z,JV)   \\ \vspace{0.1in}
		\qquad -\pi_1( X,Y,JZ,V)+2g(X,JY)g(Z,V) \}  \\ \vspace{0.1in}
		\qquad +2JZ(\nu) \{ \pi_1(X,Y,W,JV)   \\ \vspace{0.1in}
		\qquad -\pi_1( X,Y,JW,V)+2g(X,JY)g(W,V) \} =0  \,.
	\end{array} \leqno (3.2)
$$

First of all, this formula implies that $\nu$ is a constant function.

Indeed,  given a vector field V, let $Y$ be a vector field such that 
$g(V,Y)=g(JV,Y)=0$ in an open set. Putting in (3.2) $X=Z=JV$,
$W=Y$, one has:

$$
	\begin{array}{r} \vspace{0.1in}
			\ds\frac 83 (n+1)V(\nu)g(Y,Y)g(V,V)=\rho^*(R-L_3R)(JV,Y)(\nabla_V\omega)(V,Y)  \\ \vspace{0.1in}
			-\rho^*(R-L_3R)(V,Y)(\nabla_V\omega)(V,JY) \,.
	\end{array} \leqno (3.3)
$$

Therefore, $V(\nu)=0$, if $(\nabla_VJ)V=0$.

Assuming that $(\nabla_VJ)V$ does not vanish at some point, we consider an open
set where $(\nabla_VJ)V$ never vanishes and apply (3.3) to a local vector field $Y$
orthogonal to $V,\,JY,(\nabla_VJ)V,\,J(\nabla_VJ)V$. Then, we obtain again: $V(\nu)=0$.

Therefore, one has: $d\nu=0$; hence, since $M$ is connected, $\nu$ is a constant
function.

Now, the condition (3.2) turns out to be equivalent to the vanishing of the tensor
field $H$ defined by:  

$$
	\begin{array}{r} \vspace{0.1in}
			H(V,X,Y,Z,W) = \rho^*(R-L_3R)( X,Z)(\nabla_W\omega)(JY,V)  \\ \vspace{0.1in}
			- \rho^*(R-L_3R)( Y,Z)(\nabla_W\omega)(JX,V)  \\ \vspace{0.1in}
			- \rho^*(R-L_3R)(JX,Z)(\nabla_W\omega)( Y,V)  \\ \vspace{0.1in}
			+ \rho^*(R-L_3R)(JY,Z)(\nabla_W\omega)( X,V)  \\ \vspace{0.1in}
			- \rho^*(R-L_3R)( X,W)(\nabla_Z\omega)(JY,V)  \\ \vspace{0.1in}
			+ \rho^*(R-L_3R)( Y,W)(\nabla_Z\omega)(JX,V)  \\ \vspace{0.1in}
			+ \rho^*(R-L_3R)(JX,W)(\nabla_Z\omega)( Y,V)  \\ \vspace{0.1in}
			- \rho^*(R-L_3R)(JY,W)(\nabla_Z\omega)( X,V)  \,. 
	\end{array} 
$$

This implies also the vanishing of the tensor field $H'$ defined by: 

$$
	\begin{array}{r} \vspace{0.1in}
			H'(V,X,Y,Z,W) =2 (H(V,X,Y,Z,W) + H(V,Z,W,X,Y)) \\ \vspace{0.1in}
			               -H(V,Y,Z,X,W)-H(V,X,W,Y,Z) \\ \vspace{0.1in}
			               -H(V,Z,X,Y,W)-H(V,Y,W,Z,X) \,.
	\end{array} 
$$

Then, combining the conditions:

$$
	\begin{array}{c} \vspace{0.1in}
			H'(V,X,Y,Z,W) =0  \qquad H'(V,JX,JY,Z,W) =0 \\ \vspace{0.1in}
			               H'(JV,X,Y,Z,JW) =0
	\end{array} 
$$

\noindent
and using (1.9), one has:

$$
	\begin{array}{l} \vspace{0.1in}
			 \rho^*(R-L_3R)( X,Z)(\nabla_V\omega)(JY,W)  \\ \vspace{0.1in}
			\qquad - \rho^*(R-L_3R)( Y,Z)(\nabla_V\omega)(JX,W)  \\ \vspace{0.1in}
			\qquad - \rho^*(R-L_3R)(JX,Z)(\nabla_V\omega)( Y,W)  \\ \vspace{0.1in}
			\qquad + \rho^*(R-L_3R)(JY,Z)(\nabla_V\omega)( X,W) =0 \,. 
	\end{array} \leqno(3.4)
$$

This implies the K\"ahler condition, i.e. $\nabla J=0$. Indeed, if $\nabla J\ne 0$,
we consider vector fields $Y,\,V$ such that $(\nabla_VJ)Y$ never vanishes in an
open set.

Putting in (3.4) $W=Y$, $X=(\nabla_VJ)Y$, one obtains, for any $Z$, $\rho^*(R-L_3R)(JY,Z)=0$
and also $\rho^*(R-L_3R)(Y ,Z)=-\rho^*(R-L_3R)(JY,JZ)=0$. Thus, (3.4) reduces to:

$$
	\rho^*(R-L_3R)(X,Z)(\nabla_V\omega)(JY,W)-\rho^*(R-L_3R)(JX,Z)(\nabla_V\omega)(Y,W)=0  \,,
$$

or, equivalently, to:

$$
	\rho^*(R-L_3R)(X,Z)J((\nabla_V J)Y)+\rho^*(R-L_3R)(JX,Z)((\nabla_V J)Y)=0  \,.
$$

Therefore, $\rho^*(R-L_3R)$ vanishes. According to [16], this means the vanishing
of the projection $p_9(R)$. Since also $p_8(R)=p_{10}(R)=0$ (see also the Proposition 1.1),
$(M,g,J)$ turns out to be a $\mathcal R_3$-almost K\"ahler manifold with p.c.a.s.c. 
Since ${\rm dim}\, M\ge 8$, a direct application of the classification theorem in [10]
implies that $(M,g,J)$ is a K\"ahler manifold with constant holomorphic sectional
curvature. This contradicts the condition $\nabla J\ne 0$.

\vspace{0.7in}
\centerline{\large REFERENCES}

\vspace{0.2in}
\noindent
[1] V. A\textsc{postolov} -- G. G\textsc{anchev} -- S: I\textsc{vanov}: {\it Compact 
Hermitian surfaces 

\ of constant antiholomorphic sectional curvature,} Proc. Amer. Math. Soc.,

\ {\bf 125}  (1997), 3705-3714. 

\vspace{0.1in}
\noindent
[2] B. -Y. C\textsc{hen} -- K. O\textsc{giue}: {\it Some characterizations of complex space forms},  

\ Duke Math. J., {\bf 40}, 1973, 797-799.

\vspace{0.1in}\noindent
[3] M. F\textsc{alcitelli} -- A. F\textsc{arinola}: {\it Locally conformal K\"ahler manifolds 
with 

\ pointwise constant antiholomorphic sectional curvature}, Riv. Mat. Univ. 

\ Parma (4), {\bf 17}
(1991), 295-314.

\vspace{0.1in}\noindent
[4] M. F\textsc{alcitelli} -- A. F\textsc{arinola} -- S. S\textsc{alamon}: {\it Almost Hermitian 
Geometry}, 

\ Differential Geom. Appl., {\bf 4} (1994), 259-282.

\vspace{0.1in}\noindent
[5] G. G\textsc{anchev}: {\it On Bochner curvature tensor in almost Hermitian manifolds},

\ Pliska Stud. Math. Bulgar., {\bf 9} (1987), 33-43.

\vspace{0.1in}\noindent
[6] G. G\textsc{anchev} -- O. T. K\textsc{assabov}: {\it Nearly K\"ahler manifolds of constant antiho-

\ lomorphic sectional curvature}, C. R. Acad. Bulg. Sci., {\bf 35} (1982), 145-147.

\vspace{0.1in}\noindent
[7] A. G\textsc{ray}: {\it Curvature identities for Hermitian and almost Hermitian manifolds},

\ Tohoku Math. Journ. (2), {\bf 28} (1976), 601-612.

\vspace{0.1in}\noindent
[8] A. G\textsc{ray} --L. H\textsc{ervella:} {\it The sixteen classes of almost Hermitian manifolds 

\ and their linear invariants}, Ann. Mat. Pura Appl. (4), {\bf 123} (1980), 35-58.

\vspace{0.1in}\noindent
[9] O. T. K\textsc{assabov}: {\it Sur le th\'eor\`eme de F. Schur pour une vari\'et\'e presque her-

\ mitienne}, C. R. Acad. Bulg. Sci., {\bf 35} (1982), 905-908.

\vspace{0.1in}\noindent
[10] O. T. K\textsc{assabov}: {\it Almost K\"ahler manifolds of constant antiholomorphic sec-

\ \ sectional curvature}, Serdica, {\bf 9} (1983), 372-376.

\vspace{0.1in}\noindent
[11] O. T. K\textsc{assabov}: {\it $AH_3$-manifolds of constant antiholomorphic sectional cur-

\ \  vature}, Pliska Stud. Math. Bulgar., {\bf 9} (1987), 52-57.
 
\vspace{0.1in}\noindent
[12] S. K\textsc{obayashi} -- K. N\textsc{omizu}: {\it Foundations of differential geometry}, Vol. I, II

\ \ (1963), (1969) Interscience publishers.

\vspace{0.1in}\noindent
[13] T. O\textsc{guro} -- K. S\textsc{ekigawa}: {\it Non existence of almost K\"ahler structures on 

\ \ hyperbolic spaces of dimension $2n(\ge4)$}, Math. Ann., {\bf 300} (1994), 317-329.

\vspace{0.1in}\noindent
[14] G. R\textsc{izza}: {\it On almost constant type manifolds}, J. Geom., {\bf 48} (1993), 174-183.

\vspace{0.1in}\noindent
[15] S. S\textsc{alamon}: {\it Riemannian geometry and holonomy groups}, Pitman Research 

\ \ Notes Math., {\bf 201} Longman (1989).

\vspace{0.1in}\noindent
[16] F. T\textsc{ricerri} -- L. V\textsc{anhecke}: {\it Curvature tensors on almost Hermitian mani-

\ \ folds}, Trans. Amer. Math. Soc., {\bf 267} (1981), 365-398.

\vspace{0.5in}

\centerline{\it Lavoro pervenuto alla redazione il 12 luglio 1996}
\centerline{\it ed accettato per la publicazione il 1 ottobre 1997}

\vspace{0.5in}

\noindent
INDIRIZZO DEGLI AUTORI:

\noindent
M. Falcitelly -- A. Farinola -- Dipartimento di Matematica -- Universit\`a -- Via E. Orabona, 
4 - 70125 Bari, Italia

\vspace{0.1in}

\noindent
O. T. Kassabov -- Higher TTransport School (BBTY) -- ``T. Kableshkov'' -- section of Math. --
Slatina 1574 Sofia, Bulgaria

\end{document}